\documentclass[11pt,reqno]{amsart}
\usepackage[usenames]{color}
\usepackage{amsmath,verbatim,graphicx,epstopdf,enumerate}
\newcommand{\D}{\mathrm{d}}

\newcommand{\lb}{\left(}

\newcommand{\rb}{\right)}
\newcommand{\PD}{\partial}

\newcommand{\Cb}{\mathbb{C}}

\newcommand{\Rb}{\mathbb{R}}
\newcommand{\Sb}{\mathbb{S}}

\newcommand{\Beq}{\begin{equation}}
	\newcommand{\Eeq}{\end{equation}}
\newcommand{\beq}{\begin{equation*}}
	\newcommand{\eeq}{\end{equation*}}
\newcommand{\bal}{\begin{align}}
	\newcommand{\eal}{\end{align}}

\newcommand{\n}{\nabla}

\newcommand{\bp}{\begin{prob}}
	\newcommand{\ep}{\end{prob}}
\newcommand{\bpr}{\begin{proof}}
	\newcommand{\epr}{\end{proof}}

\newcommand{\bel}[1]{\begin{equation}\label{#1}}
	\newcommand{\ee}{\end{equation}}

\newtheorem{theorem}{Theorem}[section]

\numberwithin{equation}{section}

\theoremstyle{definition}

\newtheorem{remark}[theorem]{Remark}

\title[Inverse problem for quasilienar elliptic equation]{Reconstruction  for the coefficients of  a quasilinear elliptic partial differential equation}
\author[C\^{a}rstea, Nakamura and Vashisth]{C\u{a}t\u{a}lin I. C\^{a}rstea$^{*}$, Gen Nakamura$^{\dagger}$ and Manmohan Vashisth$^{\diamond}$}

\address{$^{*}$School of Mathematics, Sichuan University, Chengdu, Sichuan, 610064, P.R.China.
	\newline\indent E-mail:{\tt  catalin.carstea@gmail.com}}
\address{$^{\dagger}$Department of Mathematics, Hokkaido University, Sapporo 060-0810, Japan.
	\newline\indent E-mail:{\tt \ nakamuragenn@gmail.com}}
\address{{$^{\diamond}$ Beijing Computational Science Research Center, Beijing 100193, China.
		\newline
		\indent E-mail:{\tt\  mvashisth@csrc.ac.cn}}}

\begin{document}

	\maketitle
	\begin{abstract}
		In this paper we  consider an inverse coefficients problem for a quasilinear elliptic equation of divergence form $\nabla\cdot\vec{C}(x,\nabla u(x))=0$, in a bounded smooth domain $\Omega$. We assume that $\overrightarrow{C}(x,\vec{p})=\gamma(x)\vec{p}+\vec{b}(x)|\vec{p}|^2+\mathcal{O}(|\vec{p}|^3)$, by expanding $\overrightarrow{C}(x,\vec{p})$ around $\vec{p}=0$.  We  give a reconstruction method for  $\gamma$ and $\vec{b}$ from the Dirichlet to Neumann map defined on $\partial\Omega$. 
	\end{abstract}
	
	\indent
	${}\qquad$ \textbf{\small Keywords:}
	non-linear equation, inverse problems, reconstruction, Dirichlet to Neumann map
	
	\medskip
	${}\qquad$ \textbf{\small Mathematics Subject Classifications (2010):
	} 35J66, 65M32

	\section{Introduction and statement of the main result}
	First of all, we set up a boundary value problem for a quasilinear elliptic equation of divergence form. Let $\Omega\subset\Rb^{n}$ ( $n\geq 3$) be a bounded open set with smooth boundary $\partial\Omega$. We consider the following quasilinear elliptic boundary value problem (BVP)
	\begin{equation}\label{equation of interest}
		\begin{aligned}
			\begin{cases}
				\nabla\cdot\overrightarrow{C}(x,\nabla u(x))=0,\,\,  x\in\Omega,\\
				u(x)=\epsilon f(x),\,\,\ x\in \partial\Omega,
			\end{cases}
		\end{aligned}
	\end{equation}
	where 
	$\overrightarrow{C}(x,\nabla u(x))$ is given by 
	\begin{align}\label{definition of vector C}
		\begin{aligned}
			\overrightarrow{C}(x,\nabla u(x)):= \gamma(x)\nabla u(x)+\lvert \nabla u(x)\rvert^{2}\vec{b}(x)+\overrightarrow{R}(x,\nabla u(x))
		\end{aligned}
	\end{align}
	with $\gamma,\vec{b}\in C^{\infty}(\overline{\Omega})$ and,  
	for vector $q:=(q_{1},q_{2},\cdots,q_{n})\in\mathbb{R}^{n}$,  $\overrightarrow{R}(x,q)\in C^\infty(\overline{\Omega}\times H)$ with $H:=\{q\in {\mathbb R}^n:\,|q|\le h\}$ for a
	constant $h>0$. {Throughout this paper we assume $\gamma(x)\geq C_{1}$ for some constant $C_{1}>0$ and 
 there exists a constant $C_{2}>0$ such that 
	%This constant is defined  by the following estimate: there exists a constant $C>0$ such that if $|q|<h$
	\begin{equation}\label{estimate of R}
		\lvert\partial_q^\alpha\partial_x^\beta\overrightarrow{R}(x,q)\rvert\le C_{2}|q|^{3-|\alpha|}
	\end{equation}
holds for all $(x,q)\in \overline{\Omega}\times H$ and multi-indices $\alpha,\,\beta$ with $|\alpha|\le 3$.}
%	Throughout this paper, we will assume that $ \gamma(x)\geq C$ for some constant $C>0$.
	
	%	\begin{align}\label{definition of P(x,q)}
	%	\vec{P}(x,q):=\Big(\sum_{k,l=1}^{n}c_{kl}^{1}q_{k}q_{l},\sum_{k,l=1}^{n}c_{kl}^{2}q_{k}q_{l},\sum_{k,l=1}^{n}c_{kl}^{3}q_{k}q_{l},\cdots,\sum_{k,l=1}^{n}c_{kl}^{n}q_{k}q_{l}\Big).
	%	\end{align}
	%Here we have assumed that each $c_{kl}^i\in C^{1}(\overline\Omega)$ are real valued and

	\medskip
	Under the above setup, we have the following well-posedness result for the above (BVP) which is proved in  \cite{KN}.
	\begin{theorem}\label{well-posedness result}{\rm (\cite{KN})}
		Let $n<p<\infty$. There exist $\epsilon$ and $\delta<h/2$ such that for any $f\in W^{2-1/p,\,p}(\partial\Omega)$ satisfying $\lVert f\rVert_ {W^{2-1/p,\,p}(\partial\Omega)}<\epsilon$, the following boundary value problem 
			\begin{equation*}
		\begin{aligned}
			\begin{cases}
				\nabla\cdot\overrightarrow{C}(x,\nabla u(x))=0,\,\,  x\in\Omega,\\
				u(x)=f(x),\,\,\ x\in \partial\Omega,
			\end{cases}
		\end{aligned}
	\end{equation*}
		admits a unique solution $u$ such that $\lVert u\rVert_{W^{2,p}(\Omega)}<\delta$. Moreover, there exists $C_3>0$ independent of $f$ such that 
		\begin{equation}\label{estimate on u solution to main equation}
			\lVert u\rVert_{W^{2,p}(\Omega)}\leq C_3\lVert f\rVert_ {W^{2-1/p,p}(\partial\Omega)}.
		\end{equation}
		Here $W^{2,p}(\Omega)$ and $W^{2-1/p,\,p}(\partial\Omega)$ are the usual $L^p$-Sobolev spaces of order $2$ and $2-1/p$ in $\Omega$ and on $\partial\Omega$, respectively.  \end{theorem}
	
	Based on the well-posedness of (BVP), we define the Dirichlet to Neumann (DN in short) map $\Lambda_{\overrightarrow{C}}(\epsilon f)$ by 
	\begin{equation}\label{definition for DN map}
		\Lambda_{\overrightarrow{C}}(\epsilon f):=\nu(x)\cdot\overrightarrow{C}(x,\nabla u)|_{\partial\Omega},\,\,f\in W^{2-1/p,\,p}(\partial\Omega),
	\end{equation}
	where $u$ is solution to the (BVP) and $\nu$ is the unit normal vector of $\partial\Omega$ directed into the exterior of $\Omega$.
	
	Now we state our inverse problem.
	
	\medskip
	\noindent
	{\bf Inverse problem}: Identify $\gamma$ and $\vec{b}$ from the knowledge of  DN map $\Lambda_{\overrightarrow{C}}$.
	
	\medskip
	\noindent
	\begin{remark}
		The above (BVP) is the scalar version of displacement boundary value problem for elasticity equation and $\vec{b}$'s correspond to higher order tensors of rank 6. In material science these higher order tensors are becoming important due to the demand to investigate physical phenomena in a smaller scale (see for example \cite{Lubarda} using \cite{Ciarlet} as a guide book for nonlinear elasticity). As a consequence we need to recover these higher order tensors by solving some inverse problems. Hence we can consider our inverse problem as a toy model to reconstruct tensors up to rank 6.
	\end{remark}
	
	Concerning this inverse problem, its uniqueness is already known in (\cite{KN}). Then a next very natural question is about giving a reconstruction for identifying these $\gamma$ and $\vec{b}$. 
	
	\medskip
	Our main result in this paper is the following.
	
	\begin{theorem}\label{Main theorem}
		Knowing the DN map $\Lambda_{\overrightarrow{C}}$, we can have point-wise reconstruction  for the linear part $\gamma$ and the coefficient $\vec{b}$ of the quadratic part of $\overrightarrow{C}$. (The details of the reconstruction method will be given in the proof of this theorem see Sections 2 and 3).
	\end{theorem}
	
	Let us locate our results among the well known results on inverse problems for nonlinear scalar elliptic equations using the DN map as their measured data to identify  non-linearities or extract some information about them. The first important thing to say is that, as far as we know, the known results are about  uniqueness. The major nonlinear scalar equations which have been studied up to now are of the following forms 
	\begin{itemize}
		\item [(i)] $-\Delta u+a(x,u)=0$ (\cite{Isakov and Nachman}, \cite{Isakov and Sylvester},\cite{Sun Semilinear}),
		\item[(ii)]$ -\Delta u+b(u,\nabla u)=0$ (\cite{Isakov73}),
		\item[(iii)] $\nabla\cdot(c(x,u)\nabla u)=0$ (\cite{Sun},\cite{Sun Uhlmann}),
		\item[(v)] $\nabla\cdot(\overrightarrow{C}(x,\nabla u))=0$ (\cite{KN})
	\end{itemize}
	in $\Omega$, with some appropriate conditions on the non-linearities $a(x,u)$, $b(x,\nabla u)$, $c(x,u)$, $\overrightarrow{C}(x,\nabla u)$, and we have indicated the contributing papers in the brackets.
	It should be remarked here that the uniqueness for (ii) was even given with localized DN map. The proof in \cite{KN} had one insufficient part which can be corrected by the argument given in this paper. Our main result can be considered as a further development of \cite{KN}, giving the reconstruction of the linear part and quadratic nonlinear part of $\overrightarrow{C}(x,\nabla u)$.
	
	\medskip
	The rest of this paper is organized as follows. In Section 2 we will discuss  the $\epsilon$-expansion using which  the DN map can be linearized. The linearization of DN map is the DN map for the conductivity equation with conductivity $\gamma$. Then by the famous result \cite{Nachman_reconstruction} we reconstruct $\gamma$ and hence the remaining task is to reconstruct $\vec b$. This is done in Section 3.
	
	\section{$\epsilon$-expansion of the solution to (BVP)}
	To prove the theorem, we will use  the following $\epsilon$-expansion of solution $u$ to the (BVP)
	\begin{align}\label{epsilon expansion}
		u^{f}(x)=\epsilon u^{f}_{1}(x)+\epsilon^{2}u^{f}_{2}(x)+\mathcal{O}(\epsilon^3),   
	\end{align}
	where $u_1^f$ and $u_2^f$ are given as follows. By substituting $\n u^f(x)=\epsilon\n u_{1}^f(x)+\epsilon^{2}\n  u_{2}^f(x)+\mathcal{O}(\epsilon^{3})$ in \eqref{definition of vector C}, we get
	\begin{align*}
		\begin{aligned}
			\overrightarrow{C}(x,\nabla_{x}u^{f})
			&=\gamma(x)\nabla_{x}u^{f}(x)+\lvert\nabla u^{f}(x)\rvert^{2}\vec{b}(x)+\overrightarrow{R}(x,\nabla u^{f}(x))\\
			&=\epsilon\gamma(x)\nabla u_{1}^{f}(x)+\epsilon^{2}\lb\gamma(x)\nabla u_{2}^{f}(x)+\lvert\nabla u_{1}^{f}(x)\rvert^{2}\vec{b}(x)\rb+\mathcal{O}(\epsilon^{3}).
			%				
			%					&\ \ \ \ \ \ \ \ +
			%			\Big(\sum_{k,l=1}^{n}c^{1}_{kl}(x)\partial_{k}u^{f}(x)\partial_{l}u^{f}(x),\sum_{k,l=1}^{n}c_{kl}^{2}(x)\partial_{k}u^{f}(x)\partial_{l}u^{f}(x),\cdots,\sum_{k,l=1}^{n}c_{kl}^{n}(x)\partial_{k}u^{f}(x)\partial_{l}u^{f}(x)\Big)+o(\epsilon^{3})\qquad\quad\\
			%	&=\epsilon\gamma(x)\nabla_{x}u^{f}_{1}(x)+\epsilon^{2}\gamma(x)\nabla_{x}u^{f}_{2}(x)+\epsilon^{2}\lb\sum_{k,l=1}^{n}c_{kl}^{j}(x)\partial_{k}u^{f}_{1}(x)\partial_{l}u^{f}_{1}(x)\rb_{1\leq j\leq n}+O(\epsilon^{3}),\qquad\qquad\\
			%	&\nabla_{x}\cdot\vec{C}\lb x,\nabla_{x}u^{f}(x)\rb= \epsilon\nabla_{x}\cdot\lb \gamma(x)\nabla_{x}u^{f}_{1}\rb+\epsilon^{2}\nabla_{x}\cdot\lb\gamma(x)\nabla_{x}u^{f}_{2}\rb+\epsilon^{2}\sum_{j=1}^{n}\partial_{j}\lb\sum_{k,l=1}^{n}c_{kl}^{j}\partial_{k}u^{f}_{1}\partial_{l}u^{f}_{1}\rb+O(\epsilon^{3}).\qquad\qquad\qquad
		\end{aligned}
	\end{align*}
	Now comparing the various powers of $\epsilon$ on both sides, $u_{1}^f$ is  the solution to
	\begin{equation}\label{equation for u1}
		\begin{aligned}
			\begin{cases}
				L_{\gamma}u^{f}(x):=\nabla\cdot\lb\gamma(x)\nabla u_{1}^f(x)\rb=0,\,\,x\in\Omega,\\
				u_{1}^f(x)=f(x),\,\,x\in\partial\Omega,
			\end{cases}
		\end{aligned}
	\end{equation} 
	and $u_{2}^f$ solves 
	\begin{align}\label{equation for u2}
		\begin{aligned}
			\begin{cases}
				\nabla\cdot\lb\gamma(x)\n u_{2}^f(x)\rb+\nabla\cdot\lb\vec{b}(x)\lvert \nabla u_{1}^{f}(x)\rvert^{2}\rb=0,\,\, x\in\Omega,\\
				u_{2}^f(x)=0,\,\,x\in \partial\Omega.
			\end{cases}
		\end{aligned}
	\end{align}
{	As for the justification of the above expansion, we refer to \cite{KN}}.
	
	Next, the $\epsilon$-expansion of the DN map is 
	\begin{align}\label{NtD map in epsilon aexpansion for u i}
		\begin{aligned}
			\Lambda_{\overrightarrow{C}}(\epsilon f)\Big|_{\partial\Omega}&=\epsilon\lb\gamma(x)\partial_{\nu}u^{f}_{1}(x)\rb\Big|_{\partial\Omega} +\epsilon^{2}\lb\gamma(x)\partial_{\nu}u_{2}^f(x)+\nu(x)\cdot\vec{b}(x)\lvert\nabla_{x}u_{1}^{f}(x))\rvert^{2}\rb\Big|_{\partial\Omega}+\mathcal{O}(\epsilon^{3})\\
			&=:\epsilon g_{1}(x)+\epsilon^{2}g_{2}(x)+\mathcal{O}(\epsilon^{3}).
		\end{aligned}
	\end{align}
	Hence we can know 
	\begin{align*}
		\Lambda_{\gamma}(f):=\left.\left(\gamma(x)\partial_\nu u_{1}^{f}(x)\right)\right|_{\partial\Omega}=g_{1}(x)
	\end{align*}
	and \[\left.\left(\gamma(x)\partial_{\nu}u_{2}^f(x)+\nu(x)\cdot\vec{b}(x)\lvert\nabla_{x}u_{1}^{f}(x))\rvert^{2}\right)\right|_{\partial\Omega}=g_{2}(x).\]
	Note that $\Lambda_\gamma$ is the DN map for \eqref{equation for u1}. Also, since $W^{2-1/p,\,p}(\partial\Omega)$ is dense in the $L^2$-Sobolev space $H^{1/2}(\partial\Omega)$ of order $1/2$ on $\partial\Omega$ and the boundary value problem \ref{equation for u1} with Dirichlet data $f\in H^{1/2}(\partial\Omega)$ is well-posed in $L^2$-Sobolev space $H^1(\Omega)$ of order $1$ in $\Omega$, $\Lambda_\gamma(f)$ can be defined for $f\in H^{1/2}(\partial\Omega)$. It is well-known from the work of \cite{Nachman_reconstruction} that $\gamma$ can be reconstructed from the knowledge of $\Lambda_{\gamma}$. Once knowing $\gamma(x)$, we also know $u_{1}^{f}(x)$ in $\Omega$ for every given $f\in H^{1/2}(\partial\Omega)$.
	
	For readers' convenience, we will briefly give a summary of the reconstruction given in \cite{Nachman_reconstruction}. It consists of the following five steps:
	\begin{itemize}
		\item [{\rm Step 1.}] By the determination at the boundary, reconstruct
		$\gamma$ and $\nabla\gamma$ at $\partial\Omega$ (see for example \cite{Nakamura Tanuma}).
		\item[{\rm Step 2.}] Compute the DN map $\tilde\Lambda_q: H^{1/2}(\partial\Omega)\rightarrow H^{-1/2}(\partial\Omega)$ defined by $\tilde\Lambda_q f=\partial_\nu v\big|_{\partial\Omega}$, where $v\in H^1(\Omega)$ is the solution to boundary value problem:
		$(\Delta-q)v=0\,\,\text{in}\,\,\Omega,\,\,v\big|_{\partial\Omega}=g\in H^{1/2}(\partial\Omega)$ with $q=(\Delta\sqrt\gamma)/\sqrt\gamma$, and $H^{-1/2}(\partial\Omega)$ is the dual space of $H^{1/2}(\partial\Omega)$.
		\item[{\rm Step 3.}] For any fixed $\xi\in\mathbb{R}^n$, let $\zeta\in\mathbb{C}^n$ be such that $\zeta\cdot\zeta=0$, $(\xi+\zeta)\cdot(\xi+\zeta)=0$ and define $t(\xi,\zeta)$ by
		$$t(\xi,\zeta):=\langle(\tilde\Lambda_q-\tilde\Lambda_0)e^{-ix\cdot(\zeta+\xi)}\big|_{\partial\Omega}, (2^{-1}I+S_\zeta\tilde\Lambda_q-B_\zeta)^{-1} e^{ix\cdot\zeta}\big|_{\partial\Omega}\rangle,$$
		where $S_\zeta,\,B_\zeta$ are the traces of single layer and double layer potentials of $G_\zeta:=e^{ix\cdot\zeta}(\Delta+2i\zeta\cdot\nabla)^{-1}$ to $\partial\Omega$, respectively. Here we have denoted $\tilde\Lambda_q$ when $q=0$ by $\tilde\Lambda_0$.
		\item[{\rm Step 4.}] Compute the Fourier transform of $q$ extended by $0$ outside $\Omega$ by the inversion formula:
		$$
		\lim_{|\zeta|\rightarrow\infty} t(x,\zeta)=\int_\Omega e^{-x\cdot\xi} q(x)\,dx.
		$$
		\item[{\rm Step 5.}] Solve $(\Delta-q)z=0$ in $\Omega$, $z\big|_{\partial\Omega}=\gamma^{1/2}\big|_{\partial\Omega}$ to get
		$\gamma=z^2$.
	\end{itemize}
	
	\section{Reconstruction of $\vec{b}(x)$}
	Based on what we have obtained in the previous section, in this section we will give a reconstruction for identifying $\vec{b}(x)$. Let us start this by deriving an integral identity. Take any solution $w$ of $L_\gamma w=0$ in $\Omega$, with enough regularity, and let $\beta_{w}(x):=\gamma^{-\frac{1}{2}}(x)\chi_\Omega\vec{b}(x)\cdot\nabla w(x)$, where $\chi_\Omega$ is  the characteristic function  of $\Omega$. By multiplying \eqref{equation for u2} by $w$ and integrating over $\Omega$, we have
	\begin{equation}\label{First integral identity}
		\begin{aligned}
			&\int\limits_{}\beta_{w}(x)\gamma^{\frac{1}{2}}(x)\lvert\nabla u_{1}^{f}(x)\rvert^{2}\D x=\int\limits_{\partial\Omega} \lb\gamma(x)\partial_{\nu}u_{2}^f(x)+\nu(x)\cdot\vec{b}(x)\lvert\nabla_{x}u_{1}^{f}(x))\rvert^{2}\rb w(x)\D S_{x}\,.
		\end{aligned}
	\end{equation}
	Here and hereafter $\int\,\,dx$ denotes the integration over $\mathbb{R}^n$ and $dS_x$ denotes the standard measure on $\partial\Omega$.
	
	We will polarize \eqref{First integral identity} as follows. Consider $u_{2}(x)=u_{2}^{f+g}(x)-u_{2}^{f-g}(x)$. Then from equations \eqref{equation for u2} and \eqref{First integral identity}, we get
	\begin{align}\label{Integral identity after polarization}
		\begin{aligned}
			&4\int\limits_{}\beta_{w}(x)\gamma^{\frac{1}{2}}(x)\nabla u_{1}^{f}(x)\cdot\nabla u_{1}^{g}(x)dx=\int\limits_{\partial\Omega} \lb\gamma(x)\partial_{\nu}u_{2}(x)+4\nu(x)\cdot\vec{b}(x)\nabla_{x}u_{1}^{f}(x)\cdot\nabla u_{1}^{g}(x)\rb w(x)\D S_{x}.
		\end{aligned}
	\end{align}
	The right hand side of equation \eqref{Integral identity after polarization} is known for all $f$ and $g$.

	We can choose  $u_{1}^{f}$ and $u_{1}^{g}$ to be complex geometric optics solutions
	\begin{align}\label{expression for u1f and u1g}
		\begin{aligned}
			u_{1}^{f}(x)=v_1(x)=e^{\zeta_{1}\cdot x}\gamma^{-\frac{1}{2}}(x)\lb 1+r_{1}(x,\zeta_{1})\rb,\ \text{and}\ 
			u_{1}^{g}(x)=v_2(x)=e^{\zeta_{2}\cdot x}\gamma^{-\frac{1}{2}}(x)\lb 1+r_{2}(x,\zeta_{2})\rb
		\end{aligned}
	\end{align}
	where $r_{i}$, $i=1,2$ satisfy the equations 
	\begin{equation}\label{r-eq}
		\triangle r_i+\zeta_i\cdot\nabla r_i-q r_i=q\,\,\text{in}\,\,\mathbb{R}^n,\quad q=\frac{\triangle \gamma^{\frac{1}{2}}}{\gamma^{\frac{1}{2}}},
	\end{equation}
	and the  estimate
	\begin{align}\label{estimate for error term for u1f and u1g}
		\lVert r_{i}\rVert_{H^\sigma(\Omega)}\leq \frac{C}{\lvert \zeta_{i}\rvert},\ \text{for any } \sigma>\frac{n}{2}.
	\end{align}
	The expressions for $u_{1}^{f}$ and $u_{1}^{g}$ in \eqref{expression for u1f and u1g} and  the estimate in \eqref{estimate for error term for u1f and u1g} follow from the work of \cite{Sylvester Uhlmann}. Now 
	let  $\xi\in\Rb^{n}$ be any vector and choose $\eta,k\in \Sb^{n-1}$ such that 
	\[k\cdot\xi=k\cdot\eta=\xi\cdot\eta=0.\]   
	Using these, define $\zeta_{1},\zeta_{2}\in \Cb^{n}$ by 
	\begin{align}\label{definition of zeta i}
		\zeta_{1}:= rk-i\lb\frac{\xi}{2}+s\eta\rb,\ \ \zeta_{2}:=-
		rk-i\lb\frac{\xi}{2}-s\eta\rb,
	\end{align}
	where $r$ and $s$ are chosen such that 
	\begin{align*}
		r^{2}=\frac{\lvert\xi\rvert^{2}}{4}+s^{2}.
	\end{align*}
	With this, we have
	\begin{align*}
		\begin{aligned}
			\zeta_{i}\cdot\zeta_{i}=0,\, \zeta_1+\zeta_2=-i\xi.
		\end{aligned}
	\end{align*}

	Note that
	\begin{equation*}
		\nabla v_i= e^{\zeta_i\cdot x}\left[ \zeta_i\gamma^{-\frac{1}{2}}(1+r_i)+\nabla (\gamma^{-\frac{1}{2}})(1+r_i)+\gamma^{-\frac{1}{2}}\nabla r_i   \right],
	\end{equation*}
	so
	\begin{multline*}
		\nabla v_1\cdot\nabla v_2 =e^{-i\xi\cdot x}\left[ \gamma^{-1}\zeta_1\cdot\zeta_2+\gamma^{-\frac{1}{2}}(\zeta_1+\zeta_2)\cdot\nabla(\gamma^{-\frac{1}{2}})(1+r_1)(1+r_2)\right.\\\left.
		+|\nabla(\gamma^{-\frac{1}{2}})|^2+\gamma^{-1}(\zeta_1\cdot\nabla r_2+\zeta_2\cdot\nabla r_1)   \right]+\mathcal{O}(s^{-1})\\
		=e^{-i\xi\cdot x}\left[  \left( -\frac{1}{2}|\xi|^2 \right)\gamma^{-1} - i\gamma^{-\frac{1}{2}}\xi\cdot\nabla(\gamma^{-\frac{1}{2}}) +
		|\nabla(\gamma^{-\frac{1}{2}})|^2+\gamma^{-1}(\zeta_1\cdot\nabla r_2+\zeta_2\cdot\nabla r_1)  \right]+\mathcal{O}(s^{-1}).
	\end{multline*}
	
	Consider the the term
	\begin{equation*}
		\zeta_1\cdot\nabla r_2=(-i\xi-\zeta_2)\cdot\nabla r_2=-\xi\cdot\nabla r_2-q+\triangle r_2-q r_2=-q+\mathcal{O}(s^{-1}).
	\end{equation*}
	Then 
	\begin{equation*}
		\nabla v_1\cdot\nabla v_2 = 
		e^{-i\xi\cdot x}\left[  \left( -\frac{1}{2}|\xi|^2 \right)\gamma^{-1} + i\gamma^{-\frac{3}{2}}\xi\cdot\nabla(\gamma^{\frac{1}{2}}) +
		\gamma^{-2}|\nabla(\gamma^{\frac{1}{2}})|^2-2q \gamma^{-1}  \right]+\mathcal{O}(s^{-1}).
	\end{equation*}

	Taking the limit $s\to\infty$ in \eqref{Integral identity after polarization}, we get
	\begin{equation*}
		\int e^{-i\xi\cdot x}\beta_w \left[  \left( -\frac{1}{2}|\xi|^2 \right)\gamma^{-\frac{1}{2}} + i\gamma^{-1}\xi\cdot\nabla(\gamma^{\frac{1}{2}}) + \gamma^{-\frac{3}{2}}|\nabla(\gamma^{\frac{1}{2}})|^2-2q \gamma^{-\frac{1}{2}}  \right]dx=\text{ known.}
	\end{equation*}
	It follows that
	\begin{equation*}
		\frac{1}{2}\triangle(\gamma^{-\frac{1}{2}}\beta_w)+\nabla\cdot\left( \gamma^{-1}\nabla(\gamma^{\frac{1}{2}})\beta_w \right)+
		\left( \gamma^{-\frac{3}{2}}|\nabla(\gamma^{\frac{1}{2}})|^2-2q \gamma^{-\frac{1}{2}}\right)\beta_w = \text{ known,}
	\end{equation*}
	in $\mathbb{R}^n$, in the sense of distributions, and where we have extended $\gamma$ so that it is smooth in $\mathbb{R}^n$ and the support of $\gamma-1$ is compact. Since
	\begin{equation*}
		\triangle(\gamma^{-\frac{1}{2}}\beta_w)=\gamma^{-\frac{1}{2}}\triangle \beta_w-2\gamma^{-1}\nabla(\gamma^{\frac{1}{2}})\cdot\nabla\beta_w +\left(2\gamma^{-\frac{3}{2}}|\nabla(\gamma^{\frac{1}{2}})|^2-\gamma^{-1}\triangle (\gamma^{\frac{1}{2}})  \right)\beta_w
	\end{equation*}
	and
	\begin{equation*}
		\nabla\cdot\left( \gamma^{-1}\nabla(\gamma^{\frac{1}{2}})\beta_w \right)=
		\gamma^{-1}\nabla(\gamma^{\frac{1}{2}})\cdot\nabla\beta_w+\left(\gamma^{-1}\triangle(\gamma^{\frac{1}{2}})-2\gamma^{-\frac{3}{2}}|\nabla\gamma^{\frac{1}{2}})|^2  \right)\beta_w,
	\end{equation*}
	we can conclude that $\beta_w$ satisfies
	\begin{equation}\label{beta eq}
		\triangle \beta_w-3q\beta_w=\text{ known}\,\,\,\,\text{in}\,\,\mathbb{R}^n
	\end{equation}
	in the sense of distributions. 
	
	Next we will show that $\beta_w$ can be known. Since we do know that $\beta_w$ does exist and satisfies
	\eqref{beta eq}, we only need to show such $\beta_w$ is unique.  For this it is enough to show that
	if $f\in L^2(\mathbb{R}^n)$, with compact support, satisfies
	\begin{equation*}
		\triangle f-3qf=0\,\,\text{in}\,\,\mathbb{R}^n,
	\end{equation*}
	then $f=0$. To start proving this, note that by the interior regularity of solutions of elliptic equations, $f\in C^\infty(\mathbb{R}^n)$. Further, by recalling $f$ is compactly supported, we have $f\in C_0^\infty(\mathbb{R}^n)$. 
	
	Now by the limiting absorption principle, for any fixed $\delta>1/2$ and any given $\psi\in L^2_\delta(\mathbb{R}^n)$ there exists a unique $\phi\in L_{-\delta}^2(\mathbb{R}^n)$ such that 
	\begin{equation*}
		\triangle\varphi-3q\varphi=\psi\,\,\text{in}\,\,\mathbb{R}^n,
	\end{equation*}
	where 
	$$
	L^2_{\pm\delta}(\mathbb{R}^n):=\{\eta\in L^2_{\text{loc}}(\mathbb{R}^n): \Vert\eta\Vert_{\pm\delta}:=\big(\int_{\mathbb{R}^n}(1+|x|^2)^{\pm\delta} |\eta(x)|^2\,dx\big)^{1/2}<\infty\}
	$$
	(see Theorem 3.6 in page 413 of \cite{Yafaev} for the details). This implies
	\begin{equation*}
		\langle \psi,f\rangle=\langle\Delta\phi-3q\phi,f\rangle
		=\langle\phi,\Delta f-3q f\rangle=0.
	\end{equation*}
	Then, since $L^2_\delta(\mathbb{R}^n)$ is dense in $L^2_{\text{loc}}(\mathbb{R}^n)$, we immediately have $f=0$. Summing up we have obtained the following
	\begin{equation}\label{summary}
		\gamma^{\frac{1}{2}}\beta_{w}=\vec{b}\cdot\nabla w=\text{known}\,\,\text{in}\,\,\Omega\,\,\text{for all}\,\,w\,\,\text{solving}\,\, L_{\gamma}w=0\,\,\text{in}\,\,\Omega\,\,\text{with enough regularity}.
	\end{equation}
	
	Now let $\{w_{j}\}_{1\leq j\leq n}$ be solutions of  $L_{\gamma}w_j=0$ in $\Omega$, with enough regularity, such that $\{\n w_{j}(x)\}_{1\leq j\leq n}$ are linearly independent for a.e. every $x\in \Omega$ (see Lemma 3.1 of \cite{KN} for such $\{w_{j}\}_{1\leq j\leq n})$.
	Therefore, we have that $\vec{b}(x)\cdot \n w_{j}(x)$ is known for all $1\leq j\leq n$ and $x\in\Omega$. We will denote this known value by $F_{w_j}(x)$. Thus we have the following system of equations
	\begin{align*}
		\begin{aligned}
			\begin{bmatrix}\vspace{1mm}
				\frac{\partial w_{1}}{\partial{x_{1}}}&  \frac{\partial w_{1}}{\partial{x_{2}}}& \cdots\cdots& \frac{\partial w_{1}}{\partial{x_{n}}}\\\vspace{1mm}
				\frac{\partial w_{2}}{\partial{x_{1}}}&  \frac{\partial w_{2}}{\partial{x_{2}}}& \cdots\cdots& \frac{\partial w_{2}}{\partial{x_{n}}}\\\vspace{1mm}
				\frac{\partial w_{3}}{\partial{x_{1}}}&  \frac{\partial w_{3}}{\partial{x_{2}}}&\cdots\cdots& \frac{\partial w_{3}}{\partial{x_{n}}}\\\vspace{1mm}
				\vdots&\vdots&\vdots&\vdots\\\vspace{1mm}
				\frac{\partial w_{n}}{\partial{x_{1}}}& \frac{\partial w_{n}}{\partial{x_{2}}}& \cdots\cdots& \frac{\partial w_{n}}{\partial{x_{n}}}
			\end{bmatrix}
			\begin{bmatrix}\vspace{1.5mm}
				b_{1}(x)\\\vspace{1.5mm}
				b_{2}(x)\\\vspace{1.5mm}
				b_{3}(x)\\\vspace{1.5mm}
				\vdots\\\vspace{1.5mm}
				b_{n}(x)
			\end{bmatrix}
			=
			\begin{bmatrix}\vspace{1.5mm}
				F_{w_{1}}(x)\\\vspace{1.5mm}
				F_{w_2}(x)\\\vspace{1.5mm}
				F_{w_3}(x)\\\vspace{1.5mm}
				\vdots\\\vspace{1.5mm}
				F_{w_n}(x)
			\end{bmatrix},
			\,\,\,\,\,\, x\in\Omega.
		\end{aligned}
	\end{align*}
	Since the matrix
	\[A(x):=\lb\lb\frac{\PD w_{i}}{\PD x_{j}}\rb\rb_{1\leq i,j\leq n}\]
	\begin{comment}\begin{align*}
		\begin{aligned}
			A(x)=
			\begin{bmatrix}\vspace{1mm}
				\frac{\partial w_{1}}{\partial{x_{1}}}&  \frac{\partial w_{1}}{\partial{x_{2}}}& \cdots\cdots& \frac{\partial w_{1}}{\partial{x_{n}}}\\\vspace{1mm}
				\frac{\partial w_{2}}{\partial{x_{1}}}&  \frac{\partial w_{2}}{\partial{x_{2}}}& \cdots\cdots& \frac{\partial w_{2}}{\partial{x_{n}}}\\\vspace{1mm}
				\frac{\partial w_{3}}{\partial{x_{1}}}&  \frac{\partial w_{3}}{\partial{x_{2}}}&\cdots\cdots& \frac{\partial w_{3}}{\partial{x_{n}}}\\\vspace{1mm}
				\vdots&\vdots&\vdots&\vdots\\\vspace{1mm}
				\frac{\partial w_{n}}{\partial{x_{1}}}& \frac{\partial w_{n}}{\partial{x_{2}}}& \cdots\cdots& \frac{\partial w_{n}}{\partial{x_{n}}}
			\end{bmatrix}
		\end{aligned}
	\end{align*}
	\end{comment}
	is invertible for each $x\in\Omega$, therefore we obtain that 
	\begin{align*}
		\vec{b}(x)=A^{-1}(x)\vec{F}(x),\,\,x\in\Omega,
	\end{align*}
	where  
	\begin{align*}
		\begin{aligned}
			\vec{F}(x):=
			\begin{bmatrix}\vspace{1.5mm}
				F_{w_1}(x)\\\vspace{1.5mm}
				F_{w_2}(x)\\\vspace{1.5mm}
				F_{w_3}(x)\\\vspace{1.5mm}
				\vdots\\\vspace{1.5mm}
				F_{w_n}(x)
			\end{bmatrix}
		.
		\end{aligned}
	\end{align*}
	This gives the reconstruction  for $\vec{b}$ in $\Omega$.

	\section*{Acknowledgement} 
	The work of first author was supported by the Sichuan University.	Second author was partially supported by Grant-in-Aid for Scientific Research (15K21766, 15H05740) of the Japan Society for  the  Promotion  of  Science  doing  the  research  of  this  paper. The work of third author was supported by NSAF grant (No. U1530401).
%	For doing this research the first author, second author and third author were supported, partially supported and supported by Sichuan University, Grant-in-Aid for Scientific Research (15K21766, 15H05740) of the Japan Society for the Promotion of Science and NSAF grant (No.U1530401), respectively.

\end{document}